\documentclass[a4paper]{amsart}
\newtheorem{theorem}{Theorem}[section]
\newtheorem{lemma}[theorem]{Lemma}

\theoremstyle{definition}

\newtheorem{example}[theorem]{Example}

\newtheorem{proposition}[theorem]{Proposition}
\newtheorem{claim}[theorem]{Claim}
\usepackage[dvips]{graphicx}

\theoremstyle{remark}
\newtheorem{remark}[theorem]{Remark}

\numberwithin{equation}{section}

\begin{document}

\title[Doubly primitive knots]{Alexander polynomials of doubly primitive knots}

\author{Kazuhiro Ichihara}
\address{College of General Education, Osaka Sangyo University,
Nakagaito 3-1-1, Daito, Osaka 574-8530, Japan}
\email{ichihara@las.osaka-sandai.ac.jp}
\author{Toshio Saito}
\address{Department of Mathematics, Graduate School of Science, Osaka University,
Machikaneyama 1-1, Toyonaka, Osaka 560-0043, Japan}
\email{saito@gaia.math.wani.osaka-u.ac.jp}
\thanks{The second author is supported by the 21st Century COE program \lq\lq Towards a New Basic Science; Depth and Synthesis\rq\rq, Osaka University.}
\author{Masakazu Teragaito}
\address{Department of Mathematics and Mathematics Education, Hiroshima University,
Kagamiyama 1-1-1, Higashi-hiroshima, Japan 739-8524.}
\email{teragai@hiroshima-u.ac.jp}
\thanks{The third author is partially supported by Japan Society for the Promotion of Science, Grant-in-Aid for
Scientific Research (C), 16540071.}

\date{}



\begin{abstract}
We give a formula for Alexander polynomials of doubly primitive knots.
This also gives a practical algorithm to determine the genus of any doubly primitive knot.
\end{abstract}

\maketitle
\section{Introduction}\label{sec:intro}

Let $S$ be a standardly embedded closed orientable surface of genus two in the $3$-sphere $S^3$.
Then $S^3$ is divided into two handlebodies $H$ and $H'$ of genus two by $S$.
Let $K$ be a non-trivial knot in $S^3$, which lies on $S$.
Then $K$ is said to be \textit{doubly primitive\/} if
$K$ represents a free generator of both $\pi_1(H)$ and $\pi_1(H')$.
This notion was introduced by Berge \cite{B}, and he observed that
every doubly primitive knot admits Dehn surgery which yields a lens space.
Moreover, he suggested that if a knot admits such Dehn surgery then
the knot is doubly primitive. See also \cite[Conjecture 4.5]{G} and \cite[Problem 1.78]{Ki}.
Thus it is not too much to say that doubly primitive knots constitute one important class of knots.

Recently, Ozsv\'{a}th and Szab\'{o} \cite{OS} gave strong restrictions on the Alexander polynomials of knots
which admit lens space surgeries by using knot Floer homology.
However, as they commented, their condition is not sufficient.
For example, the Alexander polynomial of the knot $10_{132}$ satisfies their condition, but
this knot has no lens space surgery.
Also, such restrictions are obtained in \cite{K,KY,KY2,KMOS}.
In this paper, we will give an explicit formula for the Alexander polynomials
of doubly primitive knots, and recover Ozsv\'{a}th and Szab\'{o}'s condition.
Although we need to know the description of the dual knot, defined below, of a doubly primitive knot,
it is easy to calculate the Alexander polynomial by hand or computer.

To state the result, we describe the parameters of doubly primitive knots.
Let $K$ be a doubly primitive knot in $S^3$.
Suppose that $K$ lies on $S$ as above.
Then $S\cap \partial N(K)$ consists of two essential loops on $\partial N(K)$, where
$N(K)$ denotes a regular neighborhood of $K$.
The isotopy class $\gamma$ of one of these loops is called the \textit{surface slope\/}
with respect to the pair $(S,K)$. 
As Berge observed, $\gamma$-Dehn surgery on $K$ yields a lens space.
Since the surface slope depends on the position of $K$ on $S$,
it is not unique in general.
In fact, any doubly primitive knot admits at most two surface slopes.
For, a surface slope corresponds to an integer in the usual way (\cite{R}), and
any non-trivial knot admits at most two integral slopes which yield lens spaces by \cite{CGLS,M}.
Suppose that the $\gamma$-surgered manifold $K(\gamma)$ on $K$ is a lens space $L(p,q)$.
(We can assume that $0<q<p$.)
Let $K^*$ be the core of the attached solid torus of $K(\gamma)$.
We call it the \textit{associated dual knot\/} of $K$ with respect to $(S,K)$.
Let $(V_1,V_2)$ be a genus one Heegaard splitting of $L(p,q)$.
A properly embedded arc $t$ in $V_i$ is said to be \textit{trivial\/} if
$t$ is isotopic into $\partial V_i$.
Then Berge \cite{B} proved that $K^*$ can be expressed as the union of two trivial arcs in $V_1$ and $V_2$.
(Unfortunately, Berge's paper \cite{B} is unpublished, but the proof can be
found in \cite{S}.)
Furthermore, $K^*$ is isotopic to a knot $K(L(p,q);k)$ in $L(p,q)$ for some integer $k$,  $1\le k<p$,
which will be defined in Section \ref{sec:standard}.
For the doubly primitive knots constructed by Berge in \cite{B}, which are expected to give
all doubly primitive knots, there is a way to obtain such a presentation \cite{S2}.

For the triplet $(p,q,k)$, we define a Laurent polynomial.
For $i\in \{0,1,2,\dots,p-1\}$, let $\Psi(i)$ be the unique number such that
$\Psi(i)q\equiv i \pmod{p}$ and $1\le \Psi(i)\le p$.
Let $\Phi(i)=\sharp \{j\ |\ \Psi(j)<\Psi(i)\ \text{and}\ 1\le j\le k-1\}$, where $\sharp$ means the cardinality.
Then put 
\[
F(t)=\sum_{i=0}^{k-1} t^{\Phi(i)p-\Psi(i)k}
\]
and $[k]=t^{k-1}+t^{k-2}+\dots +t+1$.

\begin{theorem}\label{thm:main}
Let $K$ be a doubly primitive knot in $S^3$.
For a surface slope $\gamma$ of $K$, suppose that $K(\gamma)=L(p,q)$.
Let $K(L(p,q);k)$ be the associated dual knot in $L(p,q)$.
Then the Alexander polynomial $\Delta_K(t)$ of $K$ is equal to
$F(t)/[k]$, up to multiplication by a unit $\pm t^{n}$.
\end{theorem}

Once we have $(p,q,k)$, it is easy to calculate $F(t)$.
We will demonstrate some calculations in Section \ref{sec:example}.

Ozsv\'{a}th and Sz\'{a}bo \cite{OS} showed that
any doubly primitive knot is fibered.
Hence the degree of the Alexander polynomial of a doubly primitive knot $K$ is equal to twice the genus $g(K)$.
Thus our theorem gives a practical algorithm to determine the genus of any doubly primitive knot.

The following recovers the condition by Ozsv\'{a}th and Szab\'{o} \cite[Corollary 1.3]{OS}.

\begin{theorem}\label{thm:main2}
Let $K$ be a doubly primitive knot in $S^3$.
Then the Alexander polynomial of $K$ has the form 
\[
\Delta_K(t)=1+\sum_{i=1}^{m}(-1)^i(t^{n_i}+t^{-n_i})
\]
for some sequence $0<n_1<n_2<\dots <n_m$.
\end{theorem}

\section{Standard position}\label{sec:standard}

\begin{figure}[tb]
\includegraphics*[scale=0.6]{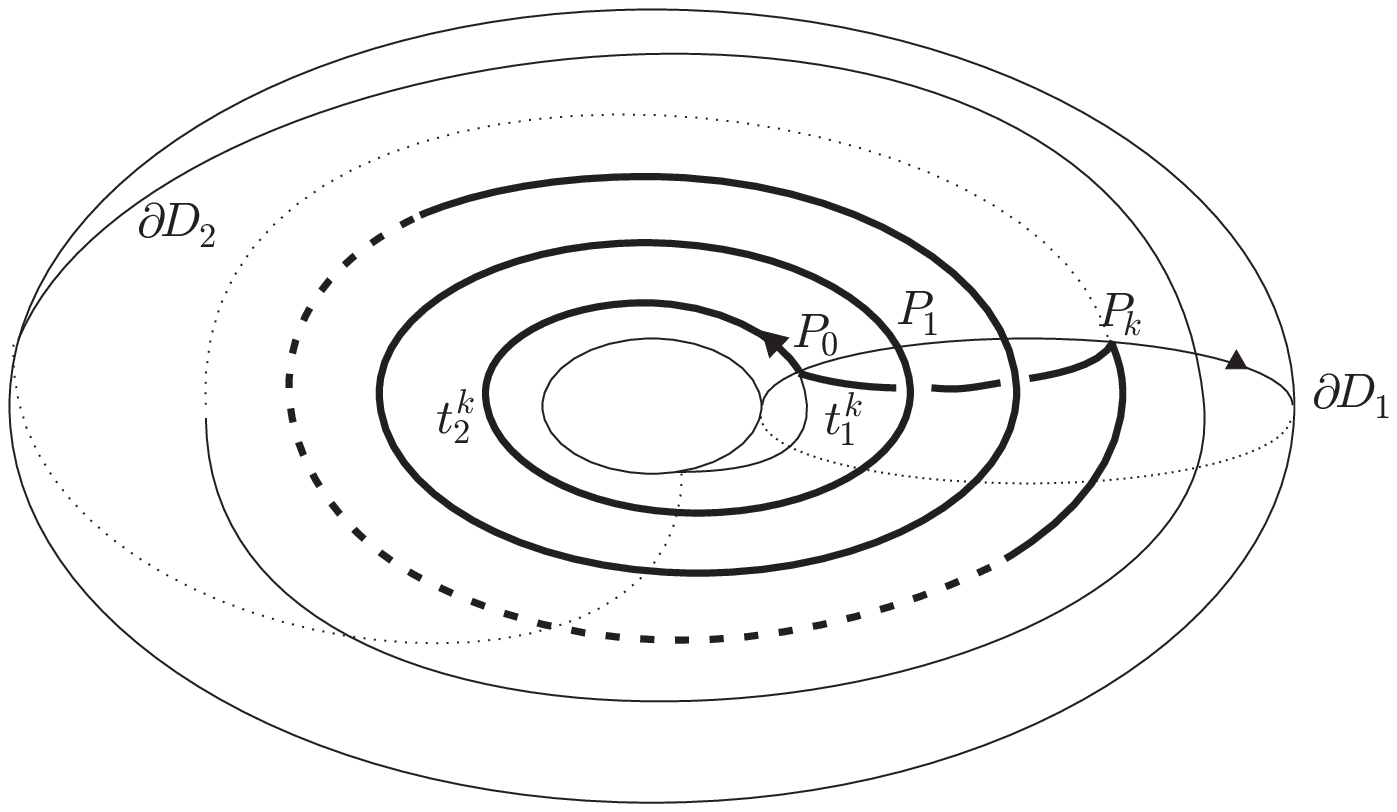}
\caption{}\label{fig:pqk}
\end{figure}

Let $K$ be a doubly primitive knot on $S$ as in Section \ref{sec:intro}.
Suppose that $K$ has the surface slope $\gamma$ with respect to $(S,K)$, and
$K(\gamma)=L(p,q)$.
Let $(V_1,V_2)$ be a genus one Heegaard splitting of $L(p,q)$.
Figure \ref{fig:pqk} shows $V_1$ with its meridian disk $D_1$ and
$\partial D_2$ on $\partial V_1$, where $D_2$ is a meridian disk of $V_2$.
We assume that $\partial D_2$ gives a $(p,q)$-curve on $\partial V_1$ with the indicated orientation.

The intersection points of $\partial D_1$ and $\partial D_2$ are labelled $P_0,P_1,\dots,P_{p-1}$ successively
along the positive direction of $\partial D_1$.
Let $k\in \{1,2,\dots,p-1\}$.
For $i=1,2$, let $t_i^k$ be a simple arc in $D_i$ joining $P_0$ to $P_k$. 
Then the knot $t_1^k\cup t_2^k$ is denoted by $K(L(p,q);k)$.
In Figure \ref{fig:pqk}, a projection of $t_2^k$ on $\partial V_1$ is illustrated.

For $i\in \{0,1,2,\dots,p-1\}$, let
$\Psi(i)$ be the unique integer such that $\Psi(i)q\equiv i\pmod{p}$ and $1\le \Psi(i)\le p$,
and let $\Phi(i)=\sharp\{j\ |\ \Psi(j)<\Psi(i)\ \text{and}\ 1\le j\le k-1\}$ as in Section \ref{sec:intro}.
Although the function $\Psi(i)$ does not depend on $k$, $\Phi(i)$ depends on $k$.
We call the sequence $\{nq \pmod{p}\}_{n=1}^p$ the \textit{basic sequence}.
Then $\Psi(i)$ indicates the position of $i$ in the basic sequence.
For convenience, let $\Psi(p)=p$ and $\Phi(p)=k-1$.
Thus $\Psi$ determines a permutation on the set $\{1,2,\dots,p\}$, since $p$ and $q$ are coprime.
We remark that $\Psi(k)=\Psi_{p,q}(k)$ and $\Phi(k)=\Phi_{p,q}(k)$ in the notation of \cite{S}.
Saito \cite[Theorem 4.5]{S} shows that 
$p\cdot \Phi(k)-k\cdot \Psi(k)=\pm 1$ or $\pm 1-p$.
In fact, this condition is necessary but not sufficient for a knot $K(L(p,q);k)$ to admit
an integral surgery yielding $S^3$.
In particular, we have:

\begin{lemma}\label{lem:saito}
$p$ and $k$ are coprime.
\end{lemma}

\section{Presentation of knot group}

For $i\in \{1,2,\dots,p\}$, let 
\begin{equation*}
E(i)=
\begin{cases}
1 \quad\text{if there is some integer $0\le j< k$ with $j\equiv iq \pmod{p}$}, \\
0 \quad\text{otherwise}.
\end{cases}
\end{equation*}
We remark that exactly $k$ terms are non-zero among $\{E(i)\}$.

\begin{lemma}[\cite{OS}]
Let $G=\pi_1(L(p,q)-K(L(p,q);k))$.  Then $G$ has a presentation $\langle X,Y\ |\ R(X,Y)\rangle$,
where $R(X,Y)=\Pi_{i=1}^p(XY^{E(i)})$ and 
the abelianizer $\mathfrak{a}: G \to \mathbb{Z}=\langle t \rangle$ sends $X$ to $t^{-k}$ and $Y$ to $t^p$.
\end{lemma}

Let $F=\langle X,Y\rangle$ be the free group generated by $\{X,Y\}$,
and let $\phi:F\to G=F/\langle R(X,Y)\rangle$ be the canonical homomorphism.
The unique extensions of $\phi$ and $\mathfrak{a}$ to the group rings are denoted by the same symbols.
Then the Alexander matrix of the above presentation of $G$ is
\[
\begin{pmatrix}
F_X(t) & F_Y(t)
\end{pmatrix}
=
\begin{pmatrix}
\mathfrak{a}\phi(\frac{\partial R(X,Y)}{\partial X}) & \mathfrak{a}\phi(\frac{\partial R(X,Y)}{\partial Y})
\end{pmatrix}.
\]
Thus the Alexander polynomial $\Delta_K(t)$ is the greatest common divisor of $F_X(t)$ and $F_Y(t)$ (cf. \cite{CF}).

\begin{lemma}\label{lem:derivative}
$\frac{\partial R(X,Y)}{\partial X}=1+\sum_{i=1}^{p-1}\prod_{j=1}^{i}XY^{E(j)}$
and $\frac{\partial R(X,Y)}{\partial Y}=\sum_{i=1}^{p}E(i)(\prod_{j=1}^{i-1}XY^{E(j)})X$.
\end{lemma}

\begin{proof}
This is a straightforward calculation.  See \cite{CF}.
\end{proof}

Suppose $E(i_j)\ne 0$ for $i_j$, where $i_1<i_2<\dots <i_k$. 
For convenience, let $i_0=1$ when $i_1>1$.
Notice that $i_k=p$.

Let $s(i)=\sum_{j=1}^{i}E(j)$ and $c(i)=-ik+ps(i)$.
Then $s(p)=k$ and $c(p)=0$.

\begin{lemma}\label{lem:degree}
$F_X(t)=\sum_{i=1}^{p}t^{c(i)}$ and
$F_Y(t)=\sum_{j=1}^k t^{c(i_j)-p}$.
\end{lemma}

\begin{proof}
By Lemma \ref{lem:derivative}, $F_X(t)=1+\sum_{i=1}^{p-1}\prod_{j=1}^{i}t^{-k+E(j)p}=1+\sum_{i=1}^{p-1}t^{-ik+ps(i)}
=\sum_{i=1}^{p}t^{-ik+ps(i)}$.
Similarly,
$\frac{\partial R(X,Y)}{\partial Y}=\sum_{j=1}^{k}(\prod_{h=1}^{i_j-1}XY^{E(h)})X$,
and hence $F_Y(t)=\sum_{j=1}^k t^{-i_jk+(j-1)p}$.
Since $c(i_j)=-i_jk+ps(i_j)=-i_jk+pj$, we have $F_Y(t)=\sum_{j=1}^k t^{c(i_j)-p}$.
\end{proof}

For positive integers $h$ and $n$,
we define $[h]^n=t^{(h-1)n}+t^{(h-2)n}+\dots+t^n+1$.
In particular, $[h]=t^{h-1}+t^{h-2}+\dots+t+1$.

\begin{lemma}\label{lem:divide}
$[p]$ divides $F_X(t)$, and $[k]$ divides $F_Y(t)$. 
\end{lemma}

\begin{proof}
Let $\zeta\ne 1$ be a $p$-th root of unity.
Then $F_X(\zeta)=\sum_{i=1}^{p}\zeta^{-ik}=0$, since $p$ and $k$ are coprime.
Hence $[p]$ divides $F_X(t)$.

Similarly, if $\xi\ne 1$ is a $k$-th root of unity, then
$F_Y(\xi)=\sum_{j=1}^{k}\xi^{(j-1)p}=0$, since $p$ and $k$ are coprime again.
Thus $[k]$ divides $F_Y(t)$.
\end{proof}

By Lemma \ref{lem:divide}, we can set $F_X(t)=[p]\ f_X(t)$ and $F_Y(t)=[k]\ f_Y(t)$.
Since $p$ and $k$ are coprime by Lemma \ref{lem:saito}, two polynomials $[p]$ and $[k]$ are also coprime.
Hence the greatest common divisor of $F_X(t)$ and $F_Y(t)$ coincides with that of $f_X(t)$ and $f_Y(t)$.


Let $d=\min\{c(i)\ |\ 1\le i\le p\}$ and $e=\min\{c(i_j)-p\ |\ 1\le j\le k\}$.
Consider the polynomial $G_X(t)=t^{-d}F_X(t)$.
Hence the lowest degree of the terms in $G_X(t)$ is zero.
Then $G_X(t)=\sum_{i=1}^{p}t^{c(i)-d}$ by Lemma \ref{lem:degree}.
A term $t^{c(i)-d}$ of $G_X(t)$ is said to be \textit{excessive\/} if $c(i)-d\ge p$.
Similarly, let $G_Y(t)=t^{-e}F_Y(t)$.
Then the lowest degree of the terms in $G_Y(t)$ is also zero, and
$G_Y(t)=\sum_{j=1}^{k}t^{c(i_j)-p-e}$.
A term $t^{c(i_j)-p-e}$ of $G_Y(t)$ is said to be \textit{excessive\/} if $c(i_j)-p-e\ge k$.
For a term $t^{c(i_j)-p-e}$ of $G_Y(t)$,
let $d(i_j)$ be the unique integer such that $0\le d(i_j)<k$ and $c(i_j)-p-e\equiv d(i_j) \pmod{k}$.
Then $(c(i_j)-p-e)-d(i_j)=m(i_j)k$ for some $m(i_j)\ge 0$.
This integer $m(i_j)$ is called the \textit{multiplicity\/} of the term.
In particular, $m(i_j)=0$ for a non-excessive term.
Since $G_Y(t)$ contains a constant term, $m(i_j)=0$ for some $j$.


The next two propositions will be proved in Section \ref{sec:calculation}.


\begin{proposition}\label{pro:gy}
$t^{-e}f_Y(t)=1+(t-1)\sum_{m(i_j)>0}t^{d(i_j)}[m(i_j)]^k$.
\end{proposition}



\begin{proposition}\label{pro:gx}
$t^{-d}f_X(t)=1+(t-1)\sum_{m(i_j)>0}t^{d(i_j)}[m(i_j)]^k$.
\end{proposition}

\begin{proof}[Proof of Theorem \ref{thm:main}]
By Propositions \ref{pro:gy} and \ref{pro:gx}, the Alexander polynomial $\Delta_K(t)=f_Y(t)=F_Y(t)/[k]$.
Hence it suffices to show that $F(t)$, defined in Section \ref{sec:intro}, coincides with $F_Y(t)$.

Let $\Psi^{-1}$ be the inverse of the permutation $\Psi$.
Then $\Psi^{-1}\{i_1,i_2,\dots,i_k\}=\{0,1,2,\dots,k-1\}$.
Since $\Phi(\Psi^{-1}(i_j))=j-1$,
$c(i_j)-p=(-i_{j}k+jp)-p=-i_jk+\Phi(\Psi^{-1}(i_j))p$.
Hence 
\begin{eqnarray*}
F(t) &=& \sum_{i=0}^{k-1}t^{\Phi(i)p-\Psi(i)k}=\sum_{j=1}^{k}t^{\Phi(\Psi^{-1}(i_j))p-\Psi(\Psi^{-1}(i_j))k} \\
     &=& \sum_{j=1}^{k}t^{c(i_j)-p}=F_Y(t)
\end{eqnarray*}
by Lemma \ref{lem:degree}.
This completes the proof of Theorem \ref{thm:main}.
\end{proof}

Similarly, we can show that $F_X(t)=\sum_{i=0}^{p-1}t^{(1-\Psi(i))k+\Phi(i)p}$.
Thus the Alexander polynomial $\Delta_K(t)$ is also equal to $F_X(t)/[p]$.
But it is simpler to use $F_Y(t)/[k]$ for a calculation, because $k<p$.

\section{Degree sequences}\label{sec:calculation}

If $1\le a<b\le p$, then the interval $[a,b]$ means the set $\{i : a\le i\le b\}\subset \{1,2,\dots,p\}$,
and $[b,a]$ means $\{i : b\le i\le p\ \text{or}\ 1\le i\le a\}$.
Moreover, if $1\le a\le p<c<a+p$, then let $[a,c]=[a,p]\cup [1,c-p]$.
Thus $[a,c]=[a,c-p]$.

\begin{lemma}\label{lem:ci}
\begin{itemize}
\item[(1)]
Let $j\in \{0,1,2,\dots,k-1\}$.
On the interval $[i_j,i_{j+1}-1]$, $c(i)=-ik+pj$, and $c(i_{j+1})=c(i_{j+1}-1)+p-k$.
\textup{(}Here, $i_{1}-1=i_k$ if $i_1=1$.\textup{)}
\item[(2)] $\{c(i)\ |\ 1\le i\le p\}=\{0,1,2,\dots,p-1\} \pmod{p}$. 
\item[(3)] Let $i_a\le i< i_{a+1}$ for $a\in \{0,1,2,\dots,k-1\}$.
Then $c(i)=c(i_j)-(i-i_j)k+(a-j)p$ for $1\le j\le k$.
\end{itemize}
\end{lemma}

\begin{proof}
(1) On the interval $[i_j,i_{j+1}-1]$, $s(i)=j$.
Hence $c(i)=-ik+pj$.
Since $s(i_{j+1})=j+1$, $c(i_{j+1})=-i_{j+1}k+p(j+1)=c(i_{j+1}-1)+p-k$.

(2) Since $c(i)\equiv -ik \pmod{p}$, the conclusion follows from the fact that $p$ and $k$ are coprime.

(3) From (1), $c(i)=-ik+pa$ and $c(i_a)=-i_ak+pa$.
Hence $c(i)=c(i_a)-(i-i_a)k$.
If $a\ge 2$, $c(i_{a})=c(i_{a-1})-(i_a-i_{a-1})k+p$, so $c(i)=c(i_{a-1})-(i-i_{a-1})k+p$.
Thus we have $c(i)=c(i_j)-(i-i_j)k+(a-j)p$ for $1\le j\le a$.

From (1), $c(i_{a+1})=c(i_{a+1}-1)+(p-k)=c(i)-(i_{a+1}-1-i)k+p-k=c(i)-(i_{a+1}-i)k+p$.
Thus $c(i)=c(i_{a+1})-(i-i_{a+1})k-p$.
For $j\ge a+2$, $c(i_j)=c(i_{a+1})-(i_j-i_{a+1})k+(j-a-1)p$.
Hence $c(i)=c(i_j)-(i-i_j)k+(a-j)p$.
\end{proof}

Thus the degree sequence $\{c(i)\}_{i=1}^p$ of $F_X(t)$ increases by $p-k$ at $i_1,i_2,\dots,i_k$,
and decreases by $k$ elsewhere.

\begin{lemma}\label{lem:d}
$\{d(i_1),d(i_2),\dots,d(i_k)\}=\{0,1,2,\dots,k-1\}$.
\end{lemma}

\begin{proof}
Assume $d(i_j)=d(i_h)$.
Then $c(i_j)\equiv c(i_h) \pmod{k}$.
Since $c(i_j)\equiv pj \pmod{k}$ and $c(i_h)\equiv ph  \pmod{k}$,
we have $pj\equiv ph \pmod{k}$.
Thus $j=h$, because $p$ and $k$ are coprime.
Hence $\{d(i_1),d(i_2),\dots,d(i_k)\}=\{0,1,2,\dots,k-1\}$.
\end{proof}

\begin{proof}[Proof of Proposition \ref{pro:gy}]
By Lemma \ref{lem:d},
\begin{eqnarray*}
G_Y(t) &=& \sum_{j=1}^k t^{c(i_j)-p-e}=\sum_{j=1}^kt^{d(i_j)+m(i_j)k}=\sum_{m(i_j)=0}t^{d(i_j)}+\sum_{m(i_j)>0}t^{d(i_j)+m(i_j)k}\\
       &=& [k]+\sum_{m(i_j)>0}\bigl(t^{d(i_j)+m(i_j)k}-t^{d(i_j)}\bigr)=[k]+\sum_{m(i_j)>0}t^{d(i_j)}\bigl(t^{m(i_j)k}-1\bigr).\\
       &=& [k]+\sum_{m(i_j)>0}t^{d(i_j)}[k](t-1)[m(i_j)]^k = [k]\Bigl(1+(t-1)\sum_{m(i_j)>0}t^{d(i_j)}[m(i_j)]^k\Bigr).\\
\end{eqnarray*}
Since $G_Y(t)=t^{-e}F_Y(t)=t^{-e}[k]f_Y(t)$, we have the conclusion.
\end{proof}


To prove Proposition \ref{pro:gx}, we need some technical lemmas.

\begin{lemma}\label{lem:impact}
Let $m(i_j)>0$ and $n\ge 0$.
\begin{itemize}
\item[(1)] $c(i_j)-d-nk\ge p$ if and only if $n\le m(i_j)-1$.
\item[(2)] $\sum_{n=0}^{m(i_j)-1}t^{c(i_j)-d-nk-p} = t^{d(i_j)}[m(i_j)]^k$.
\end{itemize}
\end{lemma}

\begin{proof}
(1) First, we claim

\begin{claim}\label{cl:min}
$d=e+k$.
\end{claim}

\begin{proof}[Proof of Claim \ref{cl:min}]
By Lemma \ref{lem:ci},
$d=c(i_\ell-1)$ for some $\ell$, and $e=c(i_\ell)-p$.
Since $c(i_\ell)=c(i_\ell-1)+(p-k)$,
$e+p=d+(p-k)$.
\end{proof}

By Claim \ref{cl:min},
$c(i_j)-d-nk-p=(c(i_j)-p-e)-(n+1)k=d(i_j)+(m(i_j)-1-n)k$.
Thus $c(i_j)-d-nk-p\ge 0$ if and only if $m(i_j)-1\ge n$.

(2) 
\begin{eqnarray*}
\sum_{n=0}^{m(i_j)-1}t^{c(i_j)-d-nk-p} &=& t^{d(i_j)+(m(i_j)-1)k}+t^{d(i_j)+(m(i_j)-2)k}+\dots +t^{d(i_j)+k}+t^{d(i_j)} \\
                                      &=& t^{d(i_j)}(t^{(m(i_j)-1)k}+t^{(m(i_j)-2)k}+\dots+t^k+1) \\
                                              &=& t^{d(i_j)}[m(i_j)]^k.
\end{eqnarray*}
\end{proof}

Let $\mathcal{E}=\{i : \text{the $i$-th term of $G_X(t)$ is excessive}\}$, and let $\mathcal{E}'$
be its complement in $\{1,2,\dots,p\}$.
Since $G_X(t)$ contains a constant term, $\mathcal{E}'\ne \emptyset$.
Then 
\begin{equation*}
G_X(t) = \sum_{i\in \mathcal{E}'}t^{c(i)-d}+\sum_{i\in \mathcal{E}}t^{c(i)-d}.
\end{equation*}
Consider a partition of $\mathcal{E}$ as follows.
Let $W(h)=\{i : hp\le c(i)-d <(h+1)p\}$ for a positive integer $h$,
and let $\mathcal{E}=W(1)\cup W(2)\cup\dots \cup W(\ell)$.

\begin{lemma}\label{lem:height}
$\ell\le k-1$.
\end{lemma}

\begin{proof}
By Lemma \ref{lem:ci}, the sequence $\{c(i)\}_{i=1}^p$ increases
only $k$ times. (If $i_1=1$, then count $c(i_k)\to c(i_1)$.)
Thus $\mathrm{deg}\,G_X(t)\le k(p-k)<kp$.
\end{proof}

\begin{lemma}\label{lem:excessive}
Let $i\in \mathcal{E}$, and
let $A=\{j : m(i_j)>0\ \text{and}\ i\in[i_j,i_j+m(i_j)-1] \}$.
Then $i\in W(h)$ if and only if $\sharp A=h$.
\end{lemma}

\begin{proof}
Choose the biggest $i_a$ with $i_a\le i$ among $\{i_1,i_2,\dots,i_k\}$.
(If $i<i_1$, then let $i_a=i_k$.)

Assume $i\in W(h)$.
Then $hp\le c(i)-d <(h+1)p$.
We distinguish two cases.

\medskip
\textit{Case \textup{1:} $a>h$.}

For $j\in \{0,1,\dots,h-1\}$,
$c(i)-d=c(i_{a-j})-e-(i-i_{a-j}+1)k+jp$ by Lemma \ref{lem:ci}(3) and Claim \ref{cl:min}.
Hence $c(i_{a-j})-e-(i-i_{a-j}+1)k+jp\ge hp$, giving $c(i_{a-j})-p-e\ge (i-i_{a-j}+1)k\ge k$.
Hence $m(i_{a-j})>0$.

Furthermore, $c(i_{a-j})-p-e=d(i_{a-j})+m(i_{a-j})k\ge (i-i_{a-j}+1)k$ gives
$(1+m(i_{a-j}))k > (i-i_{a-j}+1)k$, since $d(i_{a-j})<k$.
Thus $m(i_{a-j})>i-i_{a-j}$, so $i\in [i_{a-j},i_{a-j}+m(i_{a-j})-1]$.

Since $c(i)=c(i_{a-h})-(i-i_{a-h})k+hp$,
$c(i)-d-hp=c(i_{a-h})-d-(i-i_{a-h})k<p$.
Then $i-i_{a-h}\ge m(i_{a-h})$ by Lemma \ref{lem:impact}(1).
That is, $i\ge i_{a-h}+m(i_{a-h})$.
Hence $A=\{a-j\ :\ 0\le j\le h-1\}$.

\medskip
\textit{Case \textup{2:} $a\le h$.}

For $1\le j\le a$, we can show that $m(i_j)>0$ and $i\in [i_{j},i_{j}+m(i_{j})-1]$ exactly as in Case 1.
By Lemma \ref{lem:ci}, $c(i)=c(i_k)-ik+ap$.
If $a=h$, then $c(i_k)-d-ik=c(i)-d-hp<p$.  Thus $i\ge m(i_k)$ by Lemma \ref{lem:impact}(1).
Hence $A=\{1,2,\dots,h\}$.
Thus we suppose $a<h$.
Consider $i_{k-j}$ for $j=0,1,2,\dots,h-a-1$.
By Lemma \ref{lem:height}, $h\le \ell\le k-1$. Hence $k-(h-a-1)\ge a+2$.
Then $c(i)=c(i_{k-j})-(i-i_{k-j})k+(a-(k-j))p$ by Lemma \ref{lem:ci}(3).  
Thus $c(i_{k-j})-p-e=c(i)-d-(j+a+1)p+(i+p-i_{k-j}+1)k\ge k$, so $m(i_{k-j})>0$.

Also, $c(i_{k-j})-p-e=d(i_{k-j})+m(i_{k-j})k\ge (i+p-i_{k-j}+1)k$ gives
$m(i_{k-j})>i+p-i_{k-j}$.
Therefore $i< i_{k-j}+m(i_{k-j})-p$, which means $i\in [i_{k-j},i_{k-j}+m(i_{k-j})-1]$.

Finally, $c(i)=c(i_{k-h+a})-(i-i_{k-h+a})k+(a-(k-h+a))p=c(i_{k-h+a})-(i-i_{k-h+a})k+(h-k)p$.
Then $c(i)-d-hp=c(i_{k-h+a})-d-(i+p-i_{k-h+a})k<p$.
By Lemma \ref{lem:impact}(1), $i+p-i_{k-h+a}\ge m(i_{k-h+a})$.
That is, $i\ge i_{k-h+a}+m(i_{k-h+a})-p$.
Hence $\sharp A=h$.


\medskip

Conversely, if $\sharp A=h$, then
we can verify that $A=\{a-j : 0\le j\le h-1\}$ if $a\ge h$ or
$A=\{k-j : 0\le j\le h-a-1\}\cup \{1,2,\dots,a\}$ if $a<h$.

In the former, $c(i)=c(i_{a-h+1})-(i-i_{a-h+1})k+(h-1)p$.
Since $c(i_{a-h+1})-d-(i-i_{a-h+1})k\ge p$ by Lemma \ref{lem:impact}(1), 
$c(i)-d-hp=c(i_{a-h+1})-d-(i-i_{a-h+1})k-p\ge 0$.
Thus $c(i)-d\ge hp$.
If $c(i)-d\ge (h+1)p$, then $a-h\in A$ if $a>h$, or $k\in A$ if $a=h$, a contradiction.
Hence $i\in W(h)$.

If $a<h$, then $c(i)=c(i_{k-h+a+1})-(i-i_{k-h+a+1})k+(a-(k-h+a+1))p=c(i_{k-h+a+1})-(i-i_{k-h+a+1})k+(h-k-1)p$.
Since $c(i_{k-h+a+1})-d-(i+p-i_{k-h+a+1})k\ge p$, 
$c(i)-d-hp=c(i_{k-h+a+1})-d-(i+p-i_{k-h+a+1})k-p\ge 0$.
If $c(i)-d\ge (h+1)p$, then $k-h+a\in A$, a contradiction.
Hence $i\in W(h)$.
\end{proof}

\begin{proof}[Proof of Proposition \ref{pro:gx}]
Let $R(t)=\sum_{i\in \mathcal{E}}t^{c(i)-d}$.
Then 
\begin{equation*}
R(t)=t^p\sum_{i\in W(1)}t^{c(i)-d-p}+t^{2p}\sum_{i\in W(2)}t^{c(i)-d-2p}+\dots+t^{\ell p}\sum_{i\in W(\ell)}t^{c(i)-d-\ell p}.
\end{equation*}

Thus
\begin{eqnarray*}
G_X(t) &=& [p]+ R(t)-(\sum_{i\in W(1)}t^{c(i)-d-p}+\sum_{i\in W(2)}t^{c(i)-d-2p}+\dots+\sum_{i\in W(\ell)}t^{c(i)-d-\ell p}) \\
       &=& [p] + (t^p-1)\sum_{i\in W(1)}t^{c(i)-d-p}+(t^{2p}-1)\sum_{i\in W(2)}t^{c(i)-d-2p}+\dots \\
       & & \qquad \dots+(t^{\ell p}-1)\sum_{i\in W(\ell)}t^{c(i)-d-\ell p} \\
       &=& [p]\Bigl(1+(t-1)\sum_{i\in W(1)}t^{c(i)-d-p}+(t-1)[2]^p\sum_{i\in W(2)}t^{c(i)-d-2p}+\dots \\
       & & \qquad \dots+(t-1)[\ell]^p\sum_{i\in W(\ell)}t^{c(i)-d-\ell p}\Bigr)\\
      &=& [p] \biggl( 1+(t-1)\Bigl([1]^p\sum_{i\in W(1)}t^{c(i)-d-p}+[2]^p\sum_{i\in W(2)}t^{c(i)-d-2p}+\dots \\
        & & \qquad \dots+[\ell]^p\sum_{i\in W(\ell)}t^{c(i)-d-\ell p}\Bigr)\biggr).
\end{eqnarray*}

Let $S_h(t)=[h]^p\sum_{i\in W(h)}t^{c(i)-d-hp}$ for $1\le h\le \ell$.
Since $G_X(t)=t^{-d}F_X(t)=t^{-d}[p]f_X(t)$, we have $t^{-d}f_X(t)=G_X(t)/[p]$.
Thus it suffices to show that $\sum_{h=1}^{\ell}S_h(t)=\sum_{m(i_j)>0}t^{d(i_j)}[m(i_j)]^k$.

First, 
\begin{eqnarray*}
S_h(t) &=& (t^{(h-1)p}+t^{(h-2)p}+\dots+t^{p}+1)\sum_{i\in W(h)}t^{c(i)-d-hp} \\
       &=& \sum_{i\in W(h)}t^{c(i)-d-p}+\sum_{i\in W(h)}t^{c(i)-d-2p}+\dots+\sum_{i\in W(h)}t^{c(i)-d-hp}.
\end{eqnarray*}

Suppose $i\in W(h)$.
Let $A=\{j : m(i_j)>0 \ \text{and}\ i\in [i_j,i_j+m(i_j)-1]\}$.
By Lemma \ref{lem:excessive}, $\sharp A=h$.
Let $i_a$ be the biggest in $A$ as in the proof of Lemma \ref{lem:excessive}.
Recall that $A=\{a-j : 0\le j\le h-1\}$ if $a\ge h$ or
$A=\{k-j : 0\le j\le h-a-1\}\cup \{1,2,\dots,a\}$ if $a<h$.
Then $t^{c(i)-d-p}$ appears in $t^{d(i_a)}[m(i_a)]^k$ by Lemma \ref{lem:impact}(2).
Similarly, $t^{c(i)-d-2p}$ appears in $t^{d(i_{a-1})}[m(i_{a-1})]^k$ if $a>1$, or $t^{d(i_{k})}[m(i_{k})]^k$ if $a=1$.
Continuing this, we see that 
$t^{c(i)-d-p}, t^{c(i)-d-2p},\dots, t^{c(i)-d-hp}$
appear in $\sum_{j\in A}t^{d(i_j)}[m(i_j)]^k$ and
the correspondence is one-one.

Conversely, let $m(i_j)>0$, and
choose a term $t^{c(i_j)-d-p-nk}$ in $t^{d(i_j)}[m(i_j)]^k$.
Let $i=i_j+n$.
Then $c(i)-d=c(i_j)-d-nk\ge p$ by Lemma \ref{lem:impact}(1).
Hence $i\in W(h)$ for some $h$.
Thus the term $t^{c(i_j)-d-p-nk}$ appears in $S_h(t)$.
\end{proof}

\begin{remark}
A computer experiment suggests that $\mathcal{E}=W(1)$, that is, $\mathrm{deg}\,G_X(t)<2p$.
If this is true, then the proofs of Lemma \ref{lem:excessive} and Proposition \ref{pro:gx} would be greatly simplified.
\end{remark}

\section{Proof of Theorem \ref{thm:main2}}\label{sec:main2}

Let $\{a_1,a_2,\dots,a_h\}=\{i_j : m(i_j)>0\}$ and $a_1<a_2<\dots<a_h$.
By Propositions \ref{pro:gy} and \ref{pro:gx}, 
the Alexander polynomial of $K$ has the form
$1+(t-1)\sum_{i=1}^{h}t^{d(a_i)}[m(a_i)]^k=1+(t-1)\sum_{i=1}^{h}(t^{d(a_i)+(m(a_i)-1)k}+t^{d(a_i)+(m(a_i)-2)k}+\dots+t^{d(a_i)+k}+t^{d(a_i)})$.
Let $U_i=\{d(a_i)+(m(a_i)-j)k : 1\le j\le m(a_i)\}$ and $V_i=\{d(a_i)+(m(a_i)-j)k+1 : 1\le j\le m(a_i)\}$ for $1\le i\le h$. 
Also, let $\mathcal{U}=\cup_{i=1}^{h}U_i$ and $\mathcal{V}=\cup_{i=1}^{h}V_i$.

\begin{lemma}\label{lem:coeff}
If $i\ne j$, then $U_i\cap U_j=\emptyset$, and hence $V_i\cap V_j=\emptyset$.
\end{lemma}

\begin{proof}
Since any element of $U_i$ is congruent to $d(a_i)$ modulo $k$,
the conclusion immediately follows from Lemma \ref{lem:d}.
\end{proof}

Lemma \ref{lem:coeff} implies that any coefficient in $\Delta_K(t)$ is $\pm 1$.

\begin{lemma}
The elements of $\mathcal{U}\cup \mathcal{V}-(\mathcal{U}\cap \mathcal{V})$ have the order
$u_1<v_1<u_2<v_2<\dots<u_m<v_m$, where $u_i\in \mathcal{U}$ and $v_j\in \mathcal{V}$.
\end{lemma}

\begin{proof}
Let $u_1=\min\{d(a_i) : 1\le i\le h\}$.  Then $u_1$ is the minimal number of $\mathcal{U}$ and
$u_1\not\in \mathcal{V}$, but $u_1+1\in \mathcal{V}$.
If $u_1+1\not\in \mathcal{U}$, then let $v_1=u_1+1$.
Otherwise, $u_1+2\in\mathcal{V}$.
If $u_1+2\not\in\mathcal{U}$, then let $v_1=u_1+2$.
Continuing this process, we finally find $v_1\in \mathcal{V}-\mathcal{U}$ satisfying that $u_1<v_1$ and
there is no element of $\mathcal{U}\cup \mathcal{V}-(\mathcal{U}\cap \mathcal{V})$ between $u_1$ and $v_1$.
The same argument shows that for any $u\in \mathcal{U}-\mathcal{V}$,
the next element appears in $\mathcal{V}$.
Since $\mathcal{U}$ and $\mathcal{V}$ have the same cardinality,
the elements of $\mathcal{U}$ and $\mathcal{V}$ must alternate.
\end{proof}

\begin{proof}[Proof of Theorem \ref{thm:main2}]
Since $G_Y(t)$ contains a constant term, the corresponding $d(i_j)$ is zero.
Hence $d(a_i)\ne 0$ for $1\le i\le h$.
This means $u_1\ne 0$.
Thus the set $\mathcal{U}\cup \mathcal{V}-(\mathcal{U}\cap \mathcal{V})$ gives
all degrees of the terms in $\Delta_K(t)$ except the constant term $1$.
Thus $\Delta_K(t)=1-t^{u_1}+t^{v_1}-\dots-t^{u_m}+t^{v_m}$.
From the reciprocity of Alexander polynomials (\cite{CF}), this completes the proof of Theorem \ref{thm:main2}.
\end{proof}

\section{Examples}\label{sec:example}

\begin{example}
Let $K$ be the right-handed trefoil.
Then $5$-surgery on $K$ yields $L(5,4)$.
The associated dual knot is $K(L(5,4);2)$ as shown in \cite[Example 5.1]{S}.
Set $p=5$, $q=4$ and $k=2$.
Let us consider the basic sequence
\[
\{nq\}_{n=1}^{5}:4,3,2,1,0.
\]
Then $\Psi(i)$ is equal to the position of $i$ in this sequence, and
$\Phi(i)$ is equal to the number of terms smaller than $k$ before the term $i$ 
in the sequence.  See Table \ref{table:trefoil}.

\begin{table}[ht]
\caption{}\label{table:trefoil}
\noindent\[
\begin{array}{c|cc}
\hline
i    & 0 & 1 \\
\hline
\Psi(i) & 5  & 4  \\
\hline
\Phi(i) & 1 & 0 \\
\hline
\Phi(i)p-\Psi(i)k & -5 & -8  \\
\hline
\end{array}
\]
\end{table}

Thus $F(t)=t^{-5}+t^{-8}=t^{-8}(t^3+1)$.
Since $[k]=[2]=t+1$, $\Delta_K(t)=F(t)/[2]=t^{-8}(t^2-t+1)\doteq t^2-t+1$.

\end{example}

\begin{example}
Let $K$ be the $(-2,3,7)$-pretzel knot.
It is well known that $18$-surgery on $K$ yields $L(18,5)$.
In fact, $K$ is doubly primitive, and the associated dual knot in $L(18,5)$ is $K(L(18,5);7)$ as shown in \cite[Example 5.2]{S}.
Set $p=18, q=5$ and $k=7$.
Then the basic sequence is
\[
\{nq\}_{n=1}^{18}: 5,10,15,2,7,12,17,4,9,14,1,6,11,16,3,8,13,0.
\]
Thus we can calculate $\Psi(i)$ and $\Phi(i)$ as in Table \ref{table:pretzel}.

\begin{table}[ht]
\caption{}\label{table:pretzel}
\noindent\[
\begin{array}{c|ccccccc}
\hline
i      & 0 & 1 & 2 & 3 & 4 & 5 & 6  \\
\hline
\Psi(i)   & 18 & 11 & 4 & 15 & 8 & 1 & 12 \\
\hline
\Phi(i)   & 6 &  3 & 1 & 5 & 2 & 0 & 4 \\
\hline
\Phi(i)p-\Psi(i)k  & -18 & -23 & -10 & -15 & -20 & -7 & -12 \\
\hline
\end{array}
\]
\end{table}

Hence 
$F(t)=t^{-23}(t^5+1+t^{13}+t^8+t^3+t^{16}+t^{11})$.
Since $[k]=[7]=t^6+t^5+t^4+t^3+t^2+t+1$,
$\Delta_K(t)=F(t)/[7]\doteq 1 - t + t^3 - t^4 + t^5 -t^6 + t^7 - t^9 + t^{10}$.

\end{example}


\bibliographystyle{amsplain}

\end{document}